\documentclass{article} 
\usepackage{graphics}  
\usepackage{graphicx,graphics} 
\usepackage{amsmath}  
\usepackage{amssymb}
\usepackage{amsthm}  
\theoremstyle{plain}    
\newtheorem{thm}{Theorem}[section]
\numberwithin{equation}{section} 
\numberwithin{figure}{section} 
\theoremstyle{plain}  
\newtheorem{cor}[thm]{Corollary} 
\theoremstyle{plain}    
\newtheorem{lem}[thm]{Lemma} 
\theoremstyle{plain}    
\newtheorem*{lem*}{Lemma} 
\theoremstyle{plain}    
\newtheorem{prop}[thm]{Proposition} 
\theoremstyle{plain} 
\theoremstyle{remark}
\newtheorem{rem}[thm]{Remark}
\theoremstyle{definition}
\newtheorem{defn}[thm]{Definition}
 
\newcommand{\N}{\mathbb{N}}

\newcommand{\R}{\mathbb{R}}

\begin{document}
\cleardoublepage
\pagestyle{myheadings}
\title{On Absolutely Minimizing Lipschitz Extensions and PDE $\Delta_\infty (u) = 0$}
\author{E. Le Gruyer\thanks{Institut National des Sciences Appliqu\'{e}es, 20 Avenue des Buttes de co\"{e}smes, 
35043 Rennes cedex, France (Erwan.Le-Gruyer@insa-rennes.fr). Acknowledgement: I thank Nicoletta Tchou and  
Italo Capuzzo Dolcetta who have remotivated me on this subject.}}

\date{}
\maketitle
\markboth{E. Le Gruyer}{On Absolutely Minimizing Lipschitz Extensions }
\pagenumbering{arabic}

\begin{abstract}
We prove the existence of Absolutely Minimizing Lipschitz Extensions by a method which 
differs from those used by G. Aronsson in general metrically convex compact metric spaces 
and R. Jensen in Euclidean spaces. Assuming Jensen's hypotheses, our method yields numerical 
schemes for computing, in euclidean $\R^{n}$, the solution of viscosity of equation 
$\Delta_\infty (u) = 0$ with Dirichlet's condition. 
\end{abstract}
\section{Introduction}
To produce an optimal solution to Tietze's extension problem in general metrically convex 
compact metric space, we have introduced a class $H_{h}$ of extension schemes which solve 
the problem \cite{ARLG3}.

In this paper we first prove in section 3 that, for any continuous Dirichlet's condition $f$, 
there exists a subsequence $(H_{h(n)}(f))_{n\in \N}$ which converges to an AMLE of $f$.  
Therefore, assuming Jensen's hypotheses \cite{Jensen2}, $H_h(f)$ approachs the solution of 
viscosity of  $\Delta_\infty (u) = 0$ under Dirichlet's condition $f$ when $h$ tends to $0$.

Unfortunately it is generally hopeless to obtain a numerical approximation of this solution on, 
say, a regular grid of step $h$ by discretisation of $H_{h}$ on this grid. In fact, by such 
a discretisation we obtain an extension which is Lipschitz-optimal not for euclidean
 metric but only for the geodesic metric on the grid.
 
To overcome this difficulty we introduce in this paper an explicit scheme of extension valid 
on any finite network contained in the considered metrically convex compact metric space
and prove that the extension converges to an AMLE when the network suitably densifies 
the metric space. As a consequence, assuming Jensen's hypotheses \cite{Jensen2}, we 
obtain numerical approximations of the solution of viscosity of $\Delta_\infty (u) = 0$ 
under Dirichlet's condition. Note here that A. Oberman \cite{OBER1} has obtained  very similar 
numerical approximations in $\R^{n}$, based upon the same numerical scheme, 
proposing a proof of the convergence of the scheme based upon the $\Delta_{\infty}-$approach of the problem.

In the whole paper $(E,d)$ denotes a metrically convex compact metric space that is a 
compact lenght space with the terminology of (\cite{ARON2}, appendix). We denote 
by $\delta$ the Hausdorff metric induced by $d$ on compact non-empty subsets of $E$.\\
The second part of the paper is organized as follows.

In section 4 we prove that solutions of (\ref{intro.1}) (see below) satisfy the 
 maximum principle and, as a corollary, uniqueness of the solution.

In section 5 we prove the existence of the solution of (\ref{intro.1}) and we study 
the stability of this solution.

In section 6 we prove the existence of an AMLE as the limit of solutions of (\ref{intro.1})
for sequences $((G_n,V_n))_{n \in \N}$ which suitably densify $E$.

\begin{defn}A network on $E$ is a couple $(G,V)$ where $G \subset E$ denotes a finite 
non-empty subset of $E$ and $V$ a mapping $x\in G \rightarrow V(x) \subset G$, ($V(x)$ is 
the neighbourhood of $x$) which satisfies\\
(P1) for any $x \in G$, $x \in V(x)$;\\
(P2) for any $x$,$y \in G$, $x \in V(y)$ iff $y \in V(x)$;\\  
(P3) for any $x$,$y \in G$, there exists $ x_{1},x_{2},...,x_{n-1},x_{n} \in G$ such that 
 $x_1=x$, $x_n=y$ and $x_i \in V(x_{i+1})$ for $i =1,...,n-1$;\\ 
(P4) for any $x \in G$, any $y \in G-V(x)$ there exists $z \in V(x)$ such that $d(z,y) < d(x,y)$.\\ 
\end{defn}

To any chain such as in (P3) we associate its lenght $\sum_{i=1}^{n-1}d(x_{i},x_{i+1})$.
 We define the geodesic metric $d_g$ on $(G,V)$ by letting $d_g(x,y)$ be the infimum of 
 the lenght of chains connecting $x$ and $y$.\\ 
It follows from (P1),(P2),(P3) that $d_g$ is a metric, that $d(x,y) \leq d_g(x,y)$ 
for $x$,$y \in G$ and that $d(x,y) = d_g(x,y)$ for $x$,$y \in G$, $x \in V(y)$. 
It follows from (P2),(P3) that if $G$ has at least two elements (assumed from now on) 
then $V(x) - \{x\} \neq \ \emptyset$ for any $x \in G$. We shall denote 
$\tilde{V}(x):=V(x) - \{x\} $. Extra-condition (P4), crucial in this paper (see the end 
of theorem \ref{theoUnicite} and  theorem \ref{th.1} iii), will be used as follows:\\
for any $x\in G$, $D$ non-empty subset of $G$, $d(x,D)>0$, there exists $y\in V(x)$ such 
that $d(y,D)<d(x,D)$.

We consider the following functional equation with Dirichlet's condition :
\begin{equation}\label{intro.1}
\left\{
\begin{array}{ll}
u(x) =\mu(u;x) &  \forall x \in G-S \mbox{;}\\
u(s) =f(s)     &  \forall s \in S \mbox{.}
\end{array}
\right.
\end{equation}
Here $S$ denotes a non-empty subset of $G$, function $f$ is the Dirichlet's condition 
defined on $S$, $u$ is the numerical unknow function defined on $G$ and
\begin{equation}\label{intro.2}
 \mu(u;x) = \inf_{z \in  \tilde{V} (x)} \sup_{q \in  \tilde{V}(x)} M(u;z,q)(x) \mbox{ ;}
\end{equation}
where
\begin{equation}\label{intro.3}
 M(u;z,q)(x) := \frac{ d(x,z) u(q) + d(x,q) u(z)}{ d(x,z) + d(x,q)} \mbox{ .}
\end{equation}
\begin{rem}\label{rem001}It can be checked that
\begin{equation}\label{eq.001}
\mu(u;x) = \sup_{z \in  \tilde{V} (x)} \inf_{q \in  \tilde{V}(x)} M(u;z,q)(x) \mbox{ .}
\end{equation}
It can also be checked that
$$ J(\mu(u;x)) =\inf_{\mu \in \R} J(\mu) $$
where
$$ J(\mu(u;x)) =\sup_{z \in  \tilde{V}(x)} \frac{\mid u(z)- \mu \mid}{d(x,z)} \mbox{ .}$$
Therefore $\mu(u;x)$ is the explicit solution of the problem of minimization considered by 
A.Oberman.
\end{rem}

\section{Basics }
Let $f$ be any function from $dom(f) \subset E$ to $\R$.
We define $\kappa(f)$ by
$$
\kappa(f) := \sup_{x,y \in dom(f), x\neq y} \frac{f(x)-f(y)}{d(x,y)} \mbox{.}
$$
We call concave modulus of continuity any mapping 
$\omega : \R^{+} \rightarrow \R^{+}$ which satisfies the following:\\
(i) $\omega(0)=0$ and $\omega$ is continuous at $0$;\\
(ii) $\omega$ is increasing: 
$h_1 \leq h_2 \ \Rightarrow \omega(h_1) \leq \omega(h_2)$;\\
(iii) $\omega$ is concave.\\
We say that
 $f$ is $\Omega-$continuous iff there exists a concave modulus of continuity 
 $\omega$ such that, for any $x,y \in dom(f)$,
 \begin{equation}\label{s0.1}
 \mid f(x) - f(y)\mid \leq \omega(d(x,y)) \mbox{.}
 \end{equation}
 For such a function $f$, we denote by $\omega(f)$ the lower bound of those 
 concave moduli of continuity which satisfy (\ref{s0.1}).\\
 For any $A\subset dom(f)$, we have obviously $\omega(f\mid A)\leq \omega(f)$ 
 (symbol $\mid$ denotes restriction to).\\
 Let us restate here results of \cite{ARLG2} which are of constant use in this 
 paper.
 \begin{prop}
 Let $f$, $g$ be any two $\Omega$-continuous real-valued functions of domain $S$ and 
 let $A$ and $B$ be any two compact non-empty subsets of $S$. Then\\
 \begin{equation}\label{Introarlg2.1}
   \parallel \omega(f) -\omega(g)\parallel_{\infty,\R^{+}} \leq 2\parallel f-g \parallel_{\infty,S} \mbox{;}
   \end{equation}
  \begin{equation}\label{Introarlg2.2}
  \parallel \omega(f\mid A) -\omega(f\mid B) \parallel_{\infty,\R^{+}}\leq 4\omega(f;\delta(A,B)) \mbox{.}
  \end{equation}
  \end{prop}
 Note that (\ref{Introarlg2.1}) and (\ref{Introarlg2.2}) have been established in \cite{ARLG2} 
 for weak moduli of continuity. It is immediate that these inequalities hold for concave 
 moduli of continuity with the same constants.
 
 \begin{rem}\label{rem111}
 So, aside the obvious fact that $AM\Omega E$ (see below) are more general than AMLE, the true reason 
 why we adopt the modulus of continuity approach rather than the Lipschitz approach in this paper 
 is that there is no equivalent of (\ref{Introarlg2.1}) and (\ref{Introarlg2.2}) for Lipschitz functions. 
 \end{rem}
Now we recall Aronsson's definition of an AMLE \cite{ARON1}. Let $e$ be a Lipschitz 
extension of a Lipschitz function $f$ of compact domain.
 \begin{defn}
 We say that $e$ is an Absolutely Minimizing Lipschitz Extension of $f$ 
 if for every non empty open $D \subset E$, $D \cap dom(f) = \emptyset$ 
 we have
 $$
 \kappa(e\mid D) = \kappa(e\mid \partial D) \mbox{,}
 $$
 where $\partial D$ denotes the boundary of $D$.
 \end{defn} 
 Characterisation below has been noticed by Aronsson \cite{ARON1}.
 \begin{prop}\label{amle111}
 An extension $e$ of $f$ is AMLE iff
 for any non empty open $D \subset E$, $D \cap dom(f) = \emptyset$ we have
  \begin{equation}
   e(x) \leq \inf_{y \in \partial D}(e(y)+\kappa(e\mid \partial D)d(x,y))  \mbox{,} \forall x\in D \mbox{,}
 \end{equation}
 and
   \begin{equation}
   \sup_{y \in \partial D}(e(y)-\kappa(e\mid \partial D)d(x,y)) \leq e(x) \mbox{,} \forall x\in D \mbox{.}
 \end{equation}
 \end{prop}
 In this paper  we use a slightly more general definition (see remark \ref{rem111}). 
 Let $e$ be a continuous extension of a $\Omega$-continuous
 function $f$. 
 \begin{defn}
 We say that $e$ is an Absolutely Minimizing $\Omega$ Extension of $f$ 
 if for every non empty open subset $D$ of $E$, $D \cap dom(f) = \emptyset$ 
 we have
 $$
 \omega(e\mid D) = \omega(e\mid \partial D) \mbox{,}
 $$
 where $\partial D$ denotes the boundary of $D$.
 \end{defn}
 The analog of proposition \ref{amle111} is: 
 \begin{prop}
 An extension $e$ of $f$ is $AM\Omega E$ iff for any non empty open 
 $D \subset E$, $D \cap dom(f) = \emptyset$ we have
  \begin{equation}\label{eqHarmo.6}
   e(x) \leq \inf_{y \in \partial D}(e(y)+\omega(e\mid \partial D;d(x,y)))  \mbox{,} \forall x\in D \mbox{,}
 \end{equation}
 and
   \begin{equation}\label{eqHarmo.7}
   \sup_{y \in \partial D}(e(y)-\omega(e\mid \partial D;d(x,y))) \leq e(x) \mbox{,} \forall x\in D \mbox{.}
 \end{equation}
 \end{prop} 
 It follows from these definitions that, if $f$ is Lipschitz and $e$ is an $AM\Omega E$ of $f$, then $e$ is an 
  $AMLE$ of $f$.  
  \section{Convergence of harmonious extensions to an $AM\Omega E$}
 Let $f$ be a continuous function of closed domain $dom(f) \subset E$.\\
 For any $h>0$, for any $x \in E$ we denote by $V_{h}(x)$ 
 the closed ball of center $x$, radius $r(x)=\inf(h,d(x,dom(f)))$.\\
 By the [theorem 3.3 of \cite{ARLG3}], there exists a unique continuous extension $H_{h}(f)$  
 from $E$ to $\R$ 
 which satisfies functional equation
 \begin{equation}\label{eqHarmo}
 g(x) = \frac{1}{2}\sup_{z \in V_{h}(x)}g(z) + \frac{1}{2}\inf_{z \in V_{h}(x)}g(z) \mbox{, } \forall x \in E \mbox{.}
 \end{equation}
 Moreover, as noticed in remark 3.4 of \cite{ARLG3}, the proof of [theorem 3.3 of \cite{ARLG3}] shows that
 $$
 \omega(H_{h}(f)) = \omega(f) \mbox{.}
 $$
 Lemma \ref{lemHarmo.4} below shows that $H_{h}(f)$ is close to an $AM\Omega E$ of $f$. 
 Its proof uses the arguments of Proposition 3.9 of \cite{ARLG3}.
 \begin{lem}\label{lemHarmo.4}
 
  For any non empty open subset $D$ of $E$ , $D \cap dom(f) = \emptyset$, we have
  \begin{equation}\label{eqHarmo.1}
   H_{h}(f)(x) \leq \inf_{y \in \partial D}(H_{h}(f)(y)+\omega(H_{h}(f)\mid \partial D;d(x,y)))
    +2\omega(f;h) \mbox{,} \forall x\in D \mbox{,}
 \end{equation}
 and
   \begin{equation}\label{eqHarmo.2}
   \sup_{y \in \partial D}(H_{h}(f)(y)-\omega(H_{h}(f)\mid \partial D;d(x,y))) -2\omega(f;h)\leq H_{h}(f)(x) 
   \mbox{,} \forall x\in D \mbox{.}
 \end{equation}
 \end{lem}
 \begin{proof}
 Since the arguments are symmetric we prove only (\ref{eqHarmo.1}).\\
 Since $E$ is a compact metrically convex metric space, by theorem 3.3 \cite{ARLG3} there exists a 
 unique extension $v$ of $H_{h}(f)\mid \partial D$ in $E$ such that\\
 \begin{equation}\label{lemHarmo1}
 v(x) = \frac{1}{2}\sup_{z \in W_{h}(x)}v(z) + \frac{1}{2}\inf_{z \in W_{h}(x)}v(z) \mbox{, } \forall x \in E \mbox{.}
 \end{equation}
 where $W_{h}(x):= \{z \in E: d(x,z) \leq inf(h,d(x,\partial D)) \} \mbox{.}$\\
 Moreover\\
 $
 v(x) - v(y) \leq \omega(H_{h}(f)\mid \partial D;d(x,y)) \mbox{, }\forall x,y \in E\mbox{.}
 $\\
 In particular we have
 \begin{equation}\label{eqHarmo20}
 v(x) - v(y) \leq \omega(H_{h}(f)\mid \partial D;d(x,y)) \mbox{, }\forall y \in \partial D \mbox{, }
 \forall x \in D\mbox{.}
 \end{equation}
 Now, let us bound $\sup_{x \in D} \mid H_{h}(f)(x)-v(x) \mid$.
 By symmetry we have only to bound from above: $\Delta = \sup_{x \in D} (H_{h}(f)(x)-v(x))$. Let
 $$
 F=\{x \in \-{D}:H_{h}(f)(x)-v(x)=\Delta \} \mbox{, } M =\sup_{x \in F} H_{h}(f)(x) \mbox{,}
 $$
 and
 $$
 \tilde{F}=\{x \in F:H_{h}(f)(x)=M \} \mbox{.}
 $$
 Let $x \in \tilde{F}$ be such that 
 \begin{equation}\label{eqHarmo21}
 d(x,\partial D) =\inf_{y\in \tilde{F}}d(y,\partial D) \mbox{.}
 \end{equation}
 Let us first show that we cannot have $d(x,\partial D) >h$. Towards a contradiction 
 let us assume it is the case.\\
  Then we have $W_h(x)=\{ z \in E: d(x,z) \leq h \}$.\\
  Since $dom(f) \cap D = \emptyset$ we infer that $d(x,S) \geq d(x, \partial D) >h$, that is
  $V_h(x)=\{ z \in E: d(x,z) \leq h \}$. Therefore $V_h(x)=W_h(x)$.\\
  Now, using a similar argument to this of  theorem 3.3-uniqueness- of \cite{ARLG3}, we infer 
  that $V_h(x) \subset \tilde{F}$. It follows that  there exists $z \in \tilde{F}$ such that 
  $d(z,\partial D) < d(x,\partial D)$ which is a contradiction with definition (\ref{eqHarmo21}) of $x$. 
  
  Now, since $d(x, \partial D) \leq h$, there exists $y \in \partial D$ such that 
  $d(x,y)\leq h$ and $H_{h}(f)(y)=v(y)$.\\
  By $\Omega$-stability  of both $v$ and $H_h(f)$ we have:\\
   $H_h(f)(x)-v(x) =H_h(f)(x)- H_h(f)(y)-v(y)-v(x) \leq 2\omega(f;h)$.Therefore
   \begin{equation}\label{eqHarmo22}
   \Delta \leq 2\omega(f;h) \mbox{.}
   \end{equation}
 Now let $x \in D$ and  $y \in \partial D$. We have\\
 $H_{h}(f)(x) - H_{h}(f)(y)-\omega(H_{h}(f)\mid \partial D;d(x,y)) = A_1+A_2$\\
 where\\
 $A_1 =H_{h}(f)(x) - v(x) \mbox{,}$\\
 and\\
 $A_2= v(x)- v(y)-\omega(H_{h}(f)\mid \partial D;d(x,y)) \mbox{.}$\\
 From inequalities (\ref{eqHarmo20}) and (\ref{eqHarmo22}), we obtain 
 $A1 + A2\leq 2\omega(f;h)+0$.
 \end{proof}

 Now we prove the existence of $AM\Omega E$ in any metrically convex compact metric space.
 \begin{thm}\label{theoArmo.2}
 Let $(h(n))_{n \in \N}$ be a sequence of positive reals 
 which converges to $0$. If  sequence $(H_{h(n)}(f))_{n \in \N}$ converges uniformly 
 to a continuous extension $g$ of $f$  then $g$ is an $AM\Omega E$ of $f$.
 \end{thm}
 \begin{proof}
 By  symmetry we prove only (\ref{eqHarmo.6}).\\
 Let $D$ be non empty open subset of $E$ such that $D \cap dom(f) \neq \emptyset$. For any $\epsilon>0$, 
 there exists $N>0$ such that $ \forall n \geq N$, $\parallel H_{h(n)}(f)-g \parallel_{\infty,E} \leq \epsilon$.
 For any $x \in D$, $y \in \partial D$, $n >N$ we have\\
 $ g(x)-g(y)-\omega(g\mid\partial D;d(x,y)) \leq A_1+A_2+A_3+A_4$,
 where\\
 $A_1 = \mid g(x)-H_{h(n)}(f)(x)\mid \leq \epsilon$;\\
 $A_2 = \mid H_{h(n)}(f)(x)-H_{h(n)}(f)(y)-\omega(H_{h(n)}(f)\mid\partial D;d(x,y))\mid \leq 2\omega(f;h(n)) $;\\
 $A_3 = \mid \omega(H_{h(n)}(f)\mid\partial D;d(x,y))-\omega(g\mid\partial D;d(x,y))  \mid \leq 2\epsilon$;\\
 $A_4 = \mid H_{h(n)}(f)(y)-g(y)\mid \leq \epsilon$.\\
 The second inequality follows from Lemma \ref{lemHarmo.4} and the third one from (\ref{Introarlg2.1}).
 We obtain $ g(x)-g(y)-\omega(g\mid\partial D;d(x,y)) \leq 4\epsilon+2\omega(f;h(n))$.\\
 By letting $n \rightarrow \infty$ we have
 $$
 g(x)-g(y)-\omega(g\mid\partial D;d(x,y)) \leq  4\epsilon \mbox{,  } \forall y \in \partial D \mbox{,}
 $$
 and we obtain the stated result.
 \end{proof}
\begin{thm}\label{theoArmo.3}
For any continuous real-valued function $f$ whose domain is a compact non-empty subset of $E$, 
there exists an $AM\Omega E$ of $f$.
\end{thm}
\begin{proof}
The set $\{H_{h}(f),h>0 \}$  is equicontinuous and equibounded. Therefore, 
 by Ascoli's theorem, there exists a subsequence $(H_{h(n)}(f))_{n \in \N}$ which 
 converges uniformly to a continuous extension of $f$ which is an  $AM\Omega E$ of $f$ 
 by theorem \ref{theoArmo.2}.
\end{proof}
 \begin{rem} 
 Moreover if, for any $f$, there exists a unique $AM\Omega E$ of $f$ 
 denoted by $H(f)$ then $ \lim_{h \rightarrow 0} H_h(f) = H(f)$. 
 In this case it follows from proposition 3.9 of \cite{ARLG3} that
 $$
  \parallel H(f\mid A) - H(f\mid B)\parallel_{\infty,E} \leq 4\omega(f;\delta(A,B)) \mbox{,}
 $$
  for any non-empty compact 
 subsets $A$,$B$ of $dom(f)$. 
 \end{rem}
 \begin{rem}
 We can summarize  the difference between Jensen's proof \cite{Jensen2} of the existence of an AMLE 
 and our own proof as follows. Jensen obtains the desired AMLE as a limit of local (because solutions of PDE)
  extensions which become more and more optimally Lipschitz. We obtain the desired AMLE as a limit 
  of optimally Lipschitz extensions which become more and more local.\\
 Aronsson \cite{ARON1} (see also \cite{JUU1} and \cite{MIL1} )
 proves the existence of $AMLE$ by giving two explicit solutions:
 $$
 u = \sup\{ w : w \mbox{ AMLE } \mbox{ of } f \mbox{ from above } \} \mbox{,}
 $$
 $$
 v = \inf\{ w : w \mbox{ AMLE } \mbox{ of } f \mbox{ from below } \} \mbox{.}
 $$
 Our proof leads to less explicit but, assuming uniqueness, more constructive solutions than Aronsson's one.\\
 \end{rem}
 \begin{rem}
 Note the formal analogy between the process $u \rightarrow \Phi_h(u)$ of harmonious regularization 
 defined by 
 $$
 \Phi_h(u)(x) = \frac{1}{2}(\sup_{y\in B_h(x)}u(y) + \inf_{y\in B_h(x)}u(y))
 $$
 which deals with $PDE$ $\Delta_{\infty}u=0$ and the process $u \rightarrow \Psi_h(u)$ of harmonic 
 regularization defined by
 $$
 \Psi_h(u)(x) = \frac{\int_{B_h(x)} u(y)dy}{\int_{B_h(x)} dy}
 $$
 which deals with $PDE$ $\Delta u=0$.
 
 It is known since Gauss that any harmonic function satisfies $\Psi_h(u)=u$ for any $h>0$. 
 The analog of this result does not hold in general for the process of harmonious regularization : 
 it can be seen by numerical tests that $\Phi_h(u) \neq u$ for $u(x,y) = x^{4/3}-y^{4/3}$ even in 
 subdomains  where this function is analytic.
 
 However some functions $u$ solutions of $\Delta_{\infty}u=0$ have this property : for example 
 linear functions, $(x_1^2+x_2^2)^{1/2}$, $arctan(x_1/x_2)$, in euclidean plane, $x_1^2-x_2^2$ 
 in the plane equipped with  sup norm. 
 \end{rem}
 \begin{rem}\label{remar3.7}
 When $(E,d)$ is $(\bar{\Omega}, \parallel \mbox{ } \parallel_2)$ 
 with $\Omega$ open convex non empty subset of euclidean $\R^n$, $dom(f)=\partial\Omega$, it can be shown 
 directly (that is whithout using theorem \ref{theoArmo.3} and the equivalence between AMLE and solution 
 of viscosity of $\Delta_{\infty}u=0$) that $H_h(f)$ converges, 
 when $h$ tends to $0$, to the solution of viscosity of $\Delta_{\infty}u=0$, $u\mid \partial \Omega =f$.
 It is a consequence of Jensen's uniqueness results \cite{Jensen2} and of a Barles-Souganidis's result \cite{BARL1} :
  see Appendix.
 \end{rem}
 \begin{rem}\label{remar3.8} The results of \cite{ARLG3} and of this section hold 
 for spaces more general than compact metrically convex metric spaces. 
 They hold in compact metric spaces $(E,d)$ having the following properties:\\
 i)
 $$
   \frac{1}{2}\sup_{q\in B_h(x)}\inf_{r\in B_h(y)}d(q,r)+
     \frac{1}{2}\sup_{r\in B_h(y)}\inf_{q\in B_h(x)}d(r,q)\leq d(x,y) \mbox{,}\mbox{ for any } x,y\in E
 $$
 
 ii) For any $x\in E$ and $y \not\in B_h(x)$ there exist $z\in B_h(x)$ such that $d(y,z)<d(y,x)$.\\
 It follows that the condition of convexity on $\Omega$ assumed in remark \ref{remar3.7} can be removed.\\
 Note that conditions i) and ii) can hold in metric spaces which can be very far from metrically 
 convex metric spaces (some finite metric spaces satisfy conditions i) and ii)) : we have therefore 
 established a theorem of existence of an AMLE under weaker hypotheses than those of 
 Milman \cite{MIL1} and Juutinen \cite{JUU1} (however our result holds only for compact spaces). 
 \end{rem}
\section{Uniqueness theorem for functional equation (\ref{intro.1})}
 As usual, we first prove a maximum principle.
  \begin{thm}\label{theoUnicite} Let $f$, $g$ be any two real-valued functions both 
 of domain $S$. 
 Let $u$,$v$ be two solutions of (\ref{intro.1}) with Dirichlet's conditions $f$ and $g$ 
 respectively. Then
  \begin{equation}\label{eqUnicite1}
 \sup_{x \in G} (u(x)-v(x)) \leq \sup_{s \in S} (f(s)-g(s)) \mbox{.}
 \end{equation}
  \end{thm}
 
 \begin{proof}
 Let us set 
  $  \Delta = \sup_{x \in G} ( u(x)-v(x) ) \mbox{, }$
  $ F = \{x \in G : u(x) - v(x) = \Delta \} \mbox{, }$ 
  $ \Lambda = \sup_{x \in F} u(x) $ and
  $\tilde{F} =\{x \in F : u(x) = \Lambda \} \mbox{.}$
  
   We start our proof  by choosing $x \in \tilde{F}$ such that 
  $d(x,S) = \inf_{y \in \tilde{F}} \mbox{ }d(y,S)$.
  
 If  $d(x,S) =0$ then  equality (\ref{eqUnicite1}) is true.\\
  
 Else, let $z_{1} \in  \tilde{V} (x)$ such that
 $$
 \sup_{q \in  \tilde{V} (x)} M(v;z_{1},q)(x) = \mu(v;x) \mbox{ .}
 $$ 
 We have
 $$
 \Delta \leq \sup_{q \in  \tilde{V} (x)} M(u;z_{1},q)(x) -\sup_{q \in  \tilde{V} (x)} M(v;z_{1},q)(x) \mbox{ .}
 $$
 Let $q_{1} \in  \tilde{V} (x)$ such that
 $$
 M(u;z_{1},q_{1}) = \sup_{q \in   \tilde{V} (x)} M(u;z_{1},q) \mbox{ .}
 $$
 We have
 $$
 \Delta \leq 
  \frac{ d(x,z_{1}) (u(q_{1}) - v(q_{1}))}{ d(x,z_{1}) + d(x,q_{1})} +
  \frac{ d(x,q_{1}) (u(z_{1}) - v(z_{1}))}{ d(x,z_{1}) + d(x,q_{1})} \leq \Delta \mbox{ ,}
 $$
 from which it follows that $ z_{1},q_{1} \in F$.
 \\

 Since 
 $$
 u(x) =  \mu(u;x) 
         \leq \sup_{q \in  \tilde{V} (x)} M(u;z_{1},q)(x)
	 = M(u;z_{1},q_{1})
	 \leq u(q_{1}) \mbox{,}
 $$
 we have $q_{1} \in \tilde{F}$ and $u(q_{1})=u(x)$.
 \\
 Since  
 $$
 u(x) \leq \frac{ d(x,z_{1}) u(q_{1}) + d(x,q_{1}) u(z_{1})}{ d(x,z_{1}) + d(x,q_{1})} \mbox{ ,}
 $$
we have $u(z_{1}) = u(x)$.\\ 
 Since
\begin{equation}\label{eqUnicite2}
u(x) = u(q_{1}) = \sup_{q \in  \tilde{V} (x)} M(u;z_{1},q) 
\end{equation}
and
$$
\forall q \in \tilde{V}(x) \mbox{ , }  u(x) \geq M(u;z_{1},q)  = 
\frac{ d(x,z_{1}) u(q) + d(x,q) u(z_{1})}{ d(x,z_{1}) + d(x,q)}
$$
with $u(z_{1}) = u(x)$,\\
 we have
\begin{equation}\label{eqUnicite3}
 \forall z \in \tilde{V}(x) \mbox{ , }  u(z) \leq u(x) \mbox{.}
 \end{equation}

We finish the proof of our assertion by prooving that $u$ is constant in $V(x)$.\\
Towards a contradiction let $q \in \tilde{V}(x)$ such that $ u(q) < u(x)  $. We have

$$
u(x)  \leq   \sup_{t \in   \tilde{V}(x)} M(u;q,t) 
      = \sup_{t \in   \tilde{V}(x)} \frac{ d(x,t) u(q) + d(x,q) u(t)}{ d(x,t) + d(x,q)} \mbox{ .}
$$
Using (\ref{eqUnicite3}), we obtain
$$
u(x) \leq \sup_{t \in \tilde{V}(x)} \frac{ d(x,t) u(q) + d(x,q) u(x)}{ d(x,t) + d(x,q)} \mbox{ .}
$$
Let $c >0$ such that $u(q) = u(x) -c$, we have\\
$$
u(x) \leq \sup_{t \in \tilde{V}(x)} \frac{ d(x,t)(u(x)-c) + d(x,q) u(x)}{ d(x,t) + d(x,q)} \mbox{ .}
$$
On the other hand, we can write \\
$$
u(x) \leq \sup_{t \in \tilde{V}(x)} (u(x) - \frac{ d(x,t)c}{ d(x,t) + d(x,q)}) \mbox{ ,}
$$
that is
$$
\inf_{t \in \tilde{V}(x)} \frac{d(x,t)c}{d(x,t)+d(x,q)} \leq 0 \mbox{.}
$$
Since
$$
\inf_{t \in \tilde{V}(x)}  \frac{ d(x,t)}{ d(x,t) + d(x,q)} > 0 \mbox{ , and } c>0 \mbox{,}
$$
we obtain the desired contradiction.\\
Now, using inequality
$$
 \Delta \leq \sup_{q \in  \tilde{V}(x)} M(u;z_{1},q) -\sup_{q \in \tilde{V}(x)} M(v;z_{1},q) \mbox{ ,}
 $$ 
 and $u =$ constant in $V(x)$, we have
  $$
 \forall q \in \tilde{V}(x), \Delta \leq 
 \frac{ d(x,z_{1}) (u(q) - v(q))}{ d(x,z_{1}) + d(x,q)} +
  \frac{ d(x,q) (u(z_{1}) - v(z_{1}))}{ d(x,z_{1}) + d(x,q)} \mbox{ .}
 $$
 Therefore
 $ \Delta \leq u(q) - v(q) \leq \Delta$, $\forall q \in \tilde{V}(x)$, that is 
  $\tilde{V}(x)\subset F$.\\
 Since $u$ is constant on $V(x)$ we have $V(x) \subset \tilde{F}$.
 Since $V$ satisfies (P4), we have \\
 $ d(x,S) = \inf_{y \in \tilde{F} } d(y,S)$ 
 and
 $\inf_{q \in  \tilde{V}(x)} d(q,S)< d(x,S)$\\
  which is  a contradiction with $V(x) \subset \tilde{F}$.
 So (\ref{eqUnicite1}) is proved.
 \end{proof}
 As immediate consequences of theorem \ref{theoUnicite} we obtain Theorem 
 \ref{theoUnicite1} and corollary \ref{coroUnicite1}:
  \begin{thm}\label{theoUnicite1}
  Functional equation (\ref{intro.1}) has a unique solution.
 \end{thm}
 
 \begin{cor}\label{coroUnicite1}
 Let $u$ a solution of (\ref{intro.1}) then
 \begin{equation}\label{f.12}
 \inf_{s \in  S} f(s) \leq u(x) \leq \sup_{s \in  S} f(s) \mbox{ ,} \forall x \in G \mbox{.}
 \end{equation}
  \end{cor}


 %
 %
 %
 %
 \section{Existence and stability of solutions of (\ref{intro.1}) .}
To prove the existence of a solution of (\ref{intro.1}), we introduce a process of evolution 
$u \rightarrow \Phi(u)$ whose the stationary state $u=\Phi(u)$ is solution of (\ref{intro.1}).
 Precisely $\Phi(u) = \Psi(u;x_{N}) \circ \Psi(u;x_{N-1}) \circ ... \circ \Psi(u;x_{1})$  
 where $\{x_1,...,x_N\}$ is an enumeration of $G-S$ and $\Psi(u;x)$, $x \in G-S$ is defined as follows:\\
 \begin{equation}\label{Exist.1}
\left\{
\begin{array}{ll}
\Psi(u;x)(y) =u(y)       &  \mbox{if } y \in G-\{x\} \mbox{;}\\
\Psi(u;x)(x) =\mu(u;x)   &  \mbox{if } y = x
\end{array}
\right.
\end{equation}
 We need three lemmas useful for existence and stability.
  \begin{lem} \label{lem4.1}
   For any two scalar-valued functions $u$,$v$ of domain $G$, we have the following properties :
  \begin{equation} \label{g4.4}
   u \leq v  \Longrightarrow  \Psi(u;x) \leq \Psi(v;x) \; \mbox{, } \forall x \in G-S \mbox{ ;}
  \end{equation}
  
  \begin{equation} \label{g4.5}
   \mid \Psi(u;x)(y) - \Psi(u;x)(z) \mid  \leq \omega(u;d_{g}(y,z)) \; \mbox{, }
   \forall x \in G-S \mbox{ ,} \forall y,z \in G \mbox{ ;}
  \end{equation}
  and
  \begin{equation} \label{g4.51}
       \sup_{y \in  G}  (\Psi(u;x)(y) - \Psi(v;x)(y)) \leq \sup_{y \in  G} (u(y)-v(y)) \mbox{ . }
   \end{equation}
  
  \end{lem}
 \begin{proof}
  Let us show (\ref{g4.4}). Let $u$, $v$ scalar-valued functions of domain $G$ 
  such that $ u \leq v $. Let $x \in G-S$. We have\\
   $ \Psi (u;x)(y) = u(y) \leq v(y)= \Psi (v;x)(y)$ for $ y \in G- \{ x \} $.\\
  Since  $ \forall z, q \in \tilde{V}(x)$ we have $ M(u;z,q)(x) \leq M(v;z,q)(x) $, therefore
  $\mu(u;x) \leq \mu(v;x)$\\

  Let us show (\ref{g4.5}). \\
  It suffices to prove that  
  $ \mid \Psi (u;x)(x) - u(y) \mid  \leq \omega(u;d_{g}(x,y)) \mbox{ , } \forall y \in G $.\\
  Let $y \in G- \{ x \} $ we have two case :
  
  First case : suppose that $y \in G-\tilde{V}(x) $. \\
  Let $z_{1} \in \tilde{V}(x)$ such that $d_{g}(x,y) =  d_{g}(x,z_{1})+d_{g}(z_{1},y)$. \\
  We have
   $$
  \Psi (u;x)(x) - u(y) 
    \leq
     \sup_{q \in  \tilde{V} (x)}  \frac{ d(x,z_{1}) (u(q)-u(y)) + d(x,q) (u(z_{1})-u(y))}{ d(x,z_{1}) + d(x,q)} 
     \mbox{ .}
  $$
  By definition of $\omega(u)$ we have
   $$
  \Psi (u;x)(x) - u(y) 
    \leq \sup_{q \in  \tilde{V}  (x)} 
    ( \frac{ d(x,z_{1}) \omega(u;d_{g}(q,y)) + d(x,q) \omega(u;d_{g}(z_{1},y)))}
    { d(x,z_{1}) + d(x,q)} ) \mbox{.}
  $$
  By concavity of $\omega(u)$ we have
     $$
  \Psi (u;x)(x) - u(y) 
    \leq 
    \sup_{q \in  \tilde{V} (x)} 
   \omega(u;\frac{ d(x,z_{1})d_{g}(q,y) + d(x,q)d_{g}(z_{1},y)}{ d(x,z_{1}) + d(x,q)} ) \mbox{,}
  $$
  Since
   $$
  d(x,z_{1})d_{g}d(q,y) + d(x,q)d_{g}(z_{1},y) 
  =
  d(x,z_{1})(d_{g}(q,y) - d(x,q)) + d(x,q)d_{g}(x,y)
  $$
  by the triangle inequality, we have
   $$
  \frac{ d(x,z_{1})(d_{g}(d(q,y)) + d(x,q)(d_{g}(z_{1},y)))}{ d(x,z_{1}) + d(x,q)}
  \leq
  d_{g}(x,y) \mbox{ .}
  $$
 
 Second case. Suppose that $y \in \tilde{V}(x) $. We have
    $$
  \Psi (u;x)(x) - u(y)
    \leq 
    \sup_{q \in  \tilde{V}(x)} (M(u;y,q) - u(y))
    \leq
    \sup_{q \in  \tilde{V}(x)} \frac{ d(x,y) \omega(u;d_{g}(q,y))}{ d(x,y) + d(x,q)}
    \mbox{ .}
  $$
  By concavity we have
   $$
  \Psi (u;x)(x) - u(y)
  \leq
  \sup_{q \in  \tilde{V}(x)} \omega(u; \frac{ d_{g}(x,y)d_{g}(q,y)}{ d_{g}(x,y) + d_{g}(x,q)})\mbox{ .}
  $$
  Since
   $$
  \frac{d_{g}(q,y)}{ d_{g}(x,y) + d_{g}(x,q)} \leq 1 \mbox{ , }
  $$
  we conclude that
  $ \Psi (u;x)(x) - u(y) \leq \omega(u;d_{g}(d(x,y)) ) \mbox{ .}$ 
  The arguments to prove that
  $ u(y) -\Psi (u;x)(x)  \leq \omega(u;d_{g}(d(x,y)) ) \mbox{ }$ 
  are symmetric using (\ref{eq.001}) instead of (\ref{intro.2}). Inequality  (\ref{g4.5}) is therefore proved.\\  
Inequality (\ref{g4.51}) holds because  we have
  $\mu(u;x) - \mu(v;x) \leq  \sup_{z \in \tilde{V}(x)}(u(z)-v(z)) \mbox{.} $
  \end{proof}
  \begin{lem} \label{lem4.2}
  For any scalar-valued functions $u$,$v$ of domain $G$, we have the following properties :
  \begin{equation} \label{g4.6}
   u \leq v  \Longrightarrow  \Phi(u) \leq \Phi(v) \;\mbox{, } \forall x \in G-S \mbox{ ;}
  \end{equation}
  
  \begin{equation} \label{g4.7}
   \mid \Phi(u)(y) - \Phi(u)(z) \mid  \leq \omega(u;d_{g}(y,z))  
  \mbox{ , } \forall x \in G-S  \mbox{ , } \forall y,z \in G \mbox{ ;}
  \end{equation}
  and
  \begin{equation} \label{g4.71}
       \inf_{z \in  G} u(z)   \leq  \Phi(u)(x) \leq \sup_{z \in  G} u(z) \mbox{ , }
   \forall x \in G  \mbox{ .}
  \end{equation}
  \end{lem}
  
  \begin{proof}
  The proof is a consequence of Lemma \ref{lem4.1}.
  \end{proof}
    Now let $U_0$ be defined by
 $$
 U_0(x) = \inf_{s \in  S}( f(s) + \omega(f;d_{g}(x,s))) \mbox{ , }\forall x \in G \mbox{ .}
 $$
 Function $U_0$ looks like classical $M^c$Shane maximal Lipschitz-optimal extension 
 of $f$ on $G$. But here $U_0$ is defined with both $d$ (in $\omega(f))$ and $d_g$ (in $d_g(x,s)$). 
 Therefore we have to check that $U_0$ is an extension of $f$.
 \begin{lem}\label{lemMac1}We have
 \begin{equation}\label{Mac.1}
 U_0(s) = f(s)  \mbox{, for } s\in S \mbox{ ;}
 \end{equation}
 \begin{equation}\label{Mac.2}
 \mid U_0(x)-U_0(y) \mid \leq \omega(f;d_g(x,y)) \mbox{, for } x,y \in G \mbox{ ;}
 \end{equation}
 \begin{equation}\label{Mac.3}
 \inf_{s \in S}f(s) \leq  U_0(x)  \leq  \sup_{s \in S}f(s) \mbox{, for } x\in G \mbox{ .}
 \end{equation}
 \end{lem}
 \begin{proof}
   Let $\tilde{s} \in S$ we have
   $$
   U_0(\tilde{s})- f(\tilde{s}) \leq f(\tilde{s}) - \omega(f;d_{g}(\tilde{s},\tilde{s}))-f(\tilde{s})
    \leq 0 \mbox{ ;}
   $$
   and
   $$
      f(\tilde{s}) -U_0(\tilde{s}) = 
   \sup_{s \in  S}( f(\tilde{s})-f(s) - \omega(f;d_{g}(\tilde{s},s))) \mbox{ ,}
   $$
   therefore
      $$
      f(\tilde{s}) -U_0(\tilde{s}) \leq 
   \sup_{s \in  S}(\omega(f;d(\tilde{s},s)) - \omega(f;d_{g}(\tilde{s},s))) \mbox{ .}
   $$
   Since $ d \leq d_{g}$ we have $f(\tilde{s}) -U_0(\tilde{s}) \leq 0$
   and $f(\tilde{s})=U_0(\tilde{s})$.
   
  Let $x$,$y \in G$ we have
  $$
  U_0(x)-U_0(y) \leq \sup_{s \in  S}(f(s) + \omega(f;d_{g}(x,s))-f(s) - \omega(f;d_{g}(y,s))) \mbox{ ,}
  $$
  Therefore
   $$
  U_0(x)-U_0(y) \leq \sup_{s \in  S}(\omega(f;d_{g}(x,s))- \omega(f;d_{g}(y,s))) 
                    \leq \omega(f;d_{g}(x,y))\mbox{ .}
  $$
  
  Let $x \in G$ we have
  $$
  U_0(x) \leq \inf_{s \in S}f(s) + \sup_{s \in S}\omega(f;d_g(x,s)) \mbox{ .}
  $$
  since $\omega(f;d_g(x,s)) \leq \sup_{s_1,s_2 \in S} (f(s_1)-f(s_2))  \mbox{,}$ we have
  $$
   U_0(x) \leq \inf_{s \in S}f(s)+\sup_{s_1,s_2 \in S} (f(s_1)-f(s_2)) \leq \sup_{s \in S}f(s) \mbox{.}
  $$
  Since $\omega(f;d_g(x,s) \geq 0$, we have $U_0(x) \geq \inf_{s \in S}f(s)$.
 
 \end{proof}
 Now we are ready to prove the existence of a solution of (\ref{intro.1}).

 \begin{thm} \label{theoExist}
 Let  $(U_{n})_{n \in \N}$ the sequence defined inductively  by 
  $U_{n+1}  = \Phi(U_{n}) \mbox{, } \forall n\in\N$. This sequence 
 converges to a solution of (\ref{intro.1}) denoted by $K(f)$:

   \begin{equation} \label{g44.1}
   K(f)(s) = f(s) \mbox{ , }\forall s \in S
   \end{equation}
  \begin{equation} \label{g44.2}
    K(f)(x) = \mu(K(f);x) \mbox{ , } \forall x \in G-S \mbox{.}
  \end{equation}  
  Moreover we have
  \begin{equation}\label{g44.3}
   \mid K(f)(x) - K(f)(y) \mid \leq \omega(f;d_{g}(x,y)) \mbox{ , } \forall x \mbox{,} y \in G \mbox{ .}
  \end{equation}
    
 \end{thm}
   \begin{proof}
   
   Let us show that $(U_{n})_{n \in \N}$ is decreasing. By Lemma \ref{lem4.2}, 
   it is sufficient to prove that $U_{1} \leq U_{0}$. 
   Given an arbitrary $k \in \{1,...,,N\}$, we have
      $$
      \Psi(U_0;x_{k})(x_{k}) - U_0(x_{k}) = \mu(U_0;x_{k}) - U_0(x_{k}) \mbox{.}
      $$
   Let $s_{k} \in S$  such that
      $$
      f(s_{k}) + \omega(f;d_{g}(x_{k},s_{k})) =\inf_{s \in  S}( f(s) + \omega(f;d_{g}(x_{k},s))) \mbox{.}
      $$
      Since
      $$
      \inf_{s \in  S}( f(s) + \omega(g)(d_{g}(x_{k},y)))\leq f(s_{k}) + \omega(g)(d_{g}(y,s_{k})) \mbox{ , } 
       \forall y \in \tilde{V}(x_{k})
      $$
      we have
        $$
      \Psi(U_0;x_{k})(x_{k}) - f(s_{k}) \leq
      \inf_{z \in   \tilde{V} (x_{k})} \sup_{q \in  \tilde{V} (x_{k})}( 
      \frac{ d_{g}(x_{k},z) \omega(f;d_{g}(q,s_{k})) + d_{g}(x_{k},q) \omega(f;d_{g}(z,s_{k}))}
      { d_{g}(x_{k},z) + d_{g}(x_{k},q)} \mbox{.}
      $$
      By concavity we have
         $$
       \Psi(U_0;x_{k})(x_{k}) - f(s_{k}) \leq
       \inf_{z \in   \tilde{V} (x_{k})} \sup_{q \in  \tilde{V} (x_{k})} \omega(f; 
       \frac{ d_{g}(x_{k},z) d_{g}(q,s_{k}) + d_{g}(x_{k},q) d_{g}(z,s_{k})}
      { d_{g}(x_{k},z) + d_{g}(x_{k},q)})  \mbox{.}
        $$
      First case. If $s_{k} \not\in \tilde{V}(x_{k})$, then  
      $ \exists z' \in \tilde{V}(x_{k})$ such that $d_{g}(x_{k},s_{k}) =  d_{g}(x_{k},z')+d_{g}(z',s_{k})$.\\
      We have
         $$
       \Psi(U_0;x_{k})(x_{k}) - f(s_{k})\leq
       \sup_{q \in  \tilde{V} (x_{k})} \omega(f;
       \frac{ d_{g}(x_{k},z') d_{g}(q,s_{k}) + d_{g}(x_{k},q) d_{g}(z',s_{k})}{ d_{g}(x_{k},z') + d_{g}(x_{k},q)})
       \mbox{.}
      $$
	Since 
       $$
       d_{g}(x_{k},z') d_{g}(q,s_{k}) + d_{g}(x_{k},q) d_{g}(z',s_{k})
        =
       d_{g}(x_{k},z')( d_{g}(q,s_{k}) -d_{g}(x_{k},q) )+ d_{g}(x_{k},q)d_{g}(x_{k},s_{k})
      $$
     by the triangle inequality we have
     $$
      \Psi(U_0;x_{k})(x_{k}) - f(s_{k})
     \leq 
    \omega(f;\frac{ d_{g}(x_{k},z')d_{g}(s_{k},x_{k}) + d_{g}(x_{k},q)d_{g}(x_{k},s_{k})}
    { d_{g}(x_{k},z') + d_{g}(x_{k},q)})
     =\omega(f;d(s_{k},x_{k})) \mbox{.}
     $$
     This last inequality clearly implies
     $\Psi(U_0;x_{k})(x_{k}) \leq U_0(x_{k}) $.
     
     Second case. if $s_{k} \in \tilde{V}(x_{k})$, then
      $$
       \Psi(U_0;x_{k})(x_{k})  - f(s_{k})\leq
       \sup_{q \in  \tilde{V} (x_{k})} \omega(f;  
	\frac{ d_{g}(x_{k},s_{k}) d_{g}(q,s_{k})}{ d_{g}(x_{k},s_{k}) + d_{g}(x_{k},q)}) \mbox{.}
     $$
      Since
     $$
     \frac{ d_{g}(q,s_{k})}{ d_{g}(x_{k},s_{k}) + d_{g}(x_{k},q)} \leq 1
      $$
      we have also $\Psi(U_0;x_{k})(x_{k}) \leq U_0(x_{k}) $.
      We conclude that
      \begin{equation} \label{equa4.1}
       \forall k=1,...,N  \mbox{ , }  \Psi(U_0;x_{k}) \leq U_0 \mbox{ .}
       \end{equation}
       By Lemma \ref{lem4.1} and this last inequality, we prove inductively that\\
     
      $\forall k=1,...,N-1 
      \mbox{ , } \Psi(U_0;x_{k+1}) \mbox{,...,}\circ \Psi(U_0;x_{1}) \leq U_0  \mbox{.}$\\
      Therefore we have $U_{1} \leq U_{0}$ and we deduce from Lemma (\ref{lem4.2})
      that sequence  $(U_{n})_{ n \in \N } $ is decreasing.
      
      By Lemmas \ref{lem4.1},\ref{lem4.2},\ref{lemMac1} we prove inductively that
      $$
       U_{n} \geq \inf_{s \in  S} U_0(s) \mbox{.}
      $$
      The sequence $(U_{n})_{ n \in \N } $  is lower bounded and decreasing 
      and therefore converges to a function denoted by $K(f)$. It remains to check 
      that (\ref{g44.1}),(\ref{g44.2})  and (\ref{g44.3}) hold. 
 For any $\epsilon >0$ there exists $N \in \N$ such that
  $$
   0 \leq \sup_{y \in G-S} (U_{n}(y) - K(f)(y)) \leq \epsilon \mbox{ , } \forall n \geq N \mbox{ .}
  $$
  For $x_{k} \in G-S$ and $ n >N$ we have
  $$
  U_{n+1}(x_{k}) = \mu(\tilde{\Psi};x_{k})
  $$
  with 
  $$
  \tilde{\Psi} = 
   \Psi(U_{n};x_{k-1}) \circ \Psi(U_{n};x_{k-2}) \circ ... \circ \Psi(U_{n};x_{1}) \mbox{.}
  $$
  Since $\tilde{\Psi}(y) \leq U_{n}(y) \mbox{, } \forall y \in G-S$ we have
  $$
  0 \leq \tilde{\Psi}(y) - K(f)(y)  \leq \epsilon \mbox{.}
  $$
  We can write
  
   $$
  \mid \mu(K(f);x_{k}) - K(f)(x_{k})\mid  \leq \mid \mu(K(f);x_{k})- \mu(\tilde{\Psi};x_{k})\mid +
  \mid  U_{n+1}(x_{k})- K(f)(x_{k})\mid \mbox{ .}
  $$
  Since
  $$
 \mid \mu(K(f);x_{k})- \mu(\tilde{\Psi};x_{k})\mid \leq \sup_{y \in V(x_k)} \mid \tilde{\Psi}(y) - K(f)(y) \mid \leq \epsilon \mbox{,}
  $$
  we obtain
  $$
  \mid \mu(K(f);x_{k})-K(f)(x_{k})\mid
    \leq 2\epsilon \mbox{.}
  $$
  By lemma \ref{lemMac1}, we have $U_0(s)= f(s)$, $\forall s \in S$ and
   $U_n(s) = f(s)$, $ \forall s \in S$, and $\forall n \in \N$. Therefore 
   $K(f)$ is an extension of $f$ and we obtain the stated result.\\
  \end{proof}
  
  Now, combining theorems \ref{theoExist} and \ref{theoUnicite1}, functional equation (\ref{intro.1})
   has $K(f)$ as unique solution. 
  As a consequence of lemmas \ref{lem4.1},\ref{lem4.2},\ref{lemMac1}  and of theorems 
  \ref{theoUnicite}, \ref{theoUnicite1}
  and \ref{theoExist} we have the following properties of stability of the extension scheme $K$:
     
  \begin{thm} \label{theoStable}Let $f$, $g$ any two real-valued functions 
  both of domain $S$.
  Then
  \begin{equation}\label{stable1}
   \mid K(f)(x) - K(f)(y) \mid \leq \omega(f;d_{g}(x,y)) \mbox{ , } \forall x \mbox{,} y \in G \mbox{ ;}
  \end{equation}
  \begin{equation}\label{stable2}
   \sup_{x \in G} (K(f)(x) - K(g)(x))  \leq  \sup_{s \in S} (f(s)-g(s)) \mbox{ ;}
  \end{equation}
  for any non-empty subsets $A$, $B$ of $S$, we have 
  \begin{equation}\label{stable3}
   \sup_{x \in G} (K(f\mid A)(x) - K(f \mid B)(x))  \leq  4\omega(f;\delta_g(A,B)) \mbox{ ,}
  \end{equation}
  where $f\mid A$ and $f\mid B$ denote the restrictions of $f$ to $A$ and $B$ and 
  $\delta_g$ Haussdorff metric constructed on geodesic metric $d_g$.
  \end{thm}
 %
 %
%
%
%
%
%
%
%
%
%
 \section{Approximation of an AMLE}
 Let $f$ denote any $\Omega$-continuous real-valued 
 function whose domain is a compact non-empty subset $S$ of $E$.
 
 In this section we shall consider sequences $(G_n,V_n)_{n \in \N}$ of networks 
 having the following properties:\\
 (Q1) $\lim_{n \rightarrow \infty}r_n =0$ \\
 where $r_n :=\sup(\delta(G_n,E),\delta(S_n,S))$ and $S_n := S \cap G_n$;\\
 (Q2)  $ \lim_{n \rightarrow \infty} \rho_n =0$ where\\
  $$\rho_n =\sup_{x\in G_n}\sup_{y \in V_n(x)} d(x,y) \mbox{;}$$
 (Q3) $\lim_{n \rightarrow \infty}\parallel d_n -d \parallel =0 $\\
 where  $d_n$ denotes geodesic metric on $(G_n,V_n)$ and 
 $$
  \parallel d_n -d \parallel := \sup_{x,y\in G_n} \mid d_n(x,y) - d(x,y) \mid       \mbox{.}
 $$
 We note $b_n(x)$ the open ball of center $x \in E$, radius $r_n$, and $B_n(x)$ 
 the closed ball of center $x \in E$, radius $\rho_n$.

 Lemma \ref{lemApp1} shows that such sequences  $(G_n,V_n)_{n \in \N}$ 
 exist in any metrically convex metric space.
 \begin{lem}\label{lemApp1}
 Sequences  $(G_n,V_n)_{n \in \N}$ exist which satisfy properties (Q1),(Q2),(Q3).
 \end{lem}
 \begin{proof}
Let $(r_{n})_{n \in \N}$ and $(\rho_{n})_{n \in \N}$ be any two sequences of
  positive reals such that:
  \begin{equation}\label{amle.1}
  r_n \leq \rho_n \mbox{, } n \in \N \mbox{;}
  \end{equation}
  \begin{equation}\label{amle.2}
  \lim_{n \rightarrow \infty} \rho_{n} = 0 \mbox{;}
  \end{equation}
   \begin{equation}\label{amle.3}
  \lim_{n \rightarrow \infty} \frac{r_n}{\rho_{n}} = 0 \mbox{.}
  \end{equation}
  For any $x \in E$, $n \in \N$, let us set $b_{n}(x) := \{y \in E : d(x,y) < r_{n} \}$,
  $B_{n}(x) := \{y \in E : d(x,y) \leq \rho_{n} \}$.
  Define $(G_n,V_n)$ as follows. Since $S$ is a compact subset of $E$, $S$ is covered by balls 
  $b_n(x)$, $x \in S$. Therefore there exists $x_1,...,x_k \in S$ such that 
  $b_n(x_i), i=1,...,k$ cover $S$. Now $E- \cup_{i=1}^{k}b_n(x_i)$ is a compact subset of 
  $E$. Let $x_{k+1},...,x_m \in E-\cup_{i=1}^{k}b_n(x_i)$ such that $\cup_{i=k+1}^{m}b_n(x_i)$
  cover $E-\cup_{i=1}^{k}b_n(x_i)$. Set $G_n := \{x_1,...,x_m\}$ and, for $x \in G_n$, set 
  $V_n(x) := G_n \cap B_n(x)$. Note that, by construction, we have $\delta(S_n,S)\leq r_n$ 
  and $\delta(G_n,E)\leq r_n$. Therefore properties (Q1), (Q2) are obviously satisfied. 
  
  Now let us show that for $n \in \N$ sufficiently large we have\\
  i) $(G_n,V_n)$ is a network; \\
  ii)$\lim_{n \rightarrow \infty}\parallel d_n -d \parallel =0 $.\\
  Properties (P1) and (P2) are immediate. To prove (P4) let $x$, $y \in G_n$
  $y \not\in V_n(x)$. By metrical convexity of $E$ there exists $t \in E$ such that 
  $d(x,t)=\rho_n -r_n$ and $d(x,y) = d(y,t)+d(t,x)$. Let $z \in G_n$ such that $d(z,t)\leq r_n$.
   One has $d(x,z) \leq d(x,t)+d(t,z) \leq \rho_n -r_n +r_n =\rho_n$. Therefore $z \in V_n(x)$.
    Moreover  $d(y,z) \leq d(y,t)+d(t,z) \leq r_n +d(x,y)-\rho_n+r_n$.
    Since for $n$ sufficiently large we have $2r_n -\rho_n <0$, we infer that $d(x,z) <d(x,y)$.
    
    Now let us prove both (P3) and ii). Let $N\in \N$, $N \geq 1$, and $x$, $y\in G_n$. 
    By metrical convexity there exists elements of $E$ $y_0=x,y_{1},...,y_N=y$ such that 
    $d(y_i,y_{i+1})=d(x,y)/N$ and $d(x,y) = \sum_{i=0}^{N-1}d(y_i,y_{i+1})$.
     
    For each $i=1,...,N-1$, choose $z_i \in G_n$ such that $d(z_i,y_i) \leq r_n$. 
    We have 
    $$d(z_i,z_{i+1}) \leq d(z_i,y_{i})+d(y_i,y_{i+1})+d(y_{i+1},z_{i+1}) 
    \leq d(x,y)/N +2r_n \mbox{.}
    $$
    Now, choosing $N$ such that 
    \begin{equation}\label{eqamle3}
    \frac{d(x,y)}{N} +2r_n \leq \rho_n \mbox{,}
    \end{equation}
    we have  $z_i \in V_n(z_{i+1})$. Property (P3) is therefore proved. Moreover 
    $d_g(x,y)  \leq \sum_{i=0}^{N-1}d(z_i,z_{i+1})$. It follows that 
    $$
    d_g(x,y)-d(x,y) \leq \sum_{i=0}^{N-1}(d(z_i,z_{i+1})-d(y_i,y_{i+1})) 
    \leq \sum_{i=0}^{N-1}2r_n =2Nr_n \mbox{.}
    $$ 
    Now, for $n$ sufficiently large, one has $\rho_n > 2r_n$. Therefore (\ref{eqamle3}) is 
    satisfied by taking $N=$ the smaller integer larger than $d(x,y)/ (\rho_n-2r_n)$. It follows that 
    $$
     2N r_n \leq 2d(x,y)/((\rho_n/r_n)-2)+2r_n \mbox{.}
     $$
    
    Therefore $\parallel d_n - d\parallel \leq 2\Delta/((\rho_n/r_n)-2)+2r_n$ 
    where $\Delta $ denotes the diameter of $(E,d)$.
    Since  $\lim_{n \rightarrow \infty} r_n = 0 $
     and $\lim_{n \rightarrow \infty} r_n/\rho_{n} = 0 $, 
     lemma \ref{lemApp1} is proved.
  
  \end{proof}
 For each $n \in \N$, let us define:\\
 1) the real-valued function $f_n$ of domain $S_n$ by
 \begin{equation}
 f_n(s)= f(s) \mbox{, } x\in S_n \mbox{;}
 \end{equation}
 2) the real-valued function $W_n$ of domain $E$ by
  \begin{equation}\label{amle.7}
   W_{n}(x) = \inf_{ s \in G_{n} } ( K_{n}(f_n)(s) + \omega(f;d(s,x))) \mbox{, }  x \in E \mbox{ ;}
  \end{equation}
 where $K_n(f_n)$ ($K_n$ for short) denotes the solution of (\ref{intro.1}) for network $(G_n,V_n)$ 
 under Dirichlet's condition $f_n$.
 %
 %
 %
 %
 %
  \begin{lem}\label{lemt.1}We have
  \begin{equation}\label{eqmac1.1}
   \mid W_{n}(x) -W_{n}(y)\mid \leq \omega(f;d(x,y))\mbox{, } \forall x,y \in E\mbox{ ;}
   \end{equation}
   \begin{equation}\label{eqmac1.2}
   \forall x \in E \mbox{ , }     \inf_{ s \in dom(f) } f(s)     
    \leq  W_{n}(x)  \leq  \sup_{ s \in dom(f) }f(s) +\omega(f;r_n)    
    \mbox{ ;}
   \end{equation}
   and
    \begin{equation}\label{eqmac1.3}
    K_{n}(s) - \sup_{ t \in G_{n} } \omega(f;d_{n}(t,s)-d(t,s))    
    \leq  W_{n}(s)  \leq  K_{n}(f_n)(s) \mbox{,  } \forall s \in G_{n}
    \mbox{ ;}
   \end{equation}
   from which we infer
   \begin{equation}\label{eqmac1.4}
    K_{n}(s) -  \omega(f;\parallel d_{n}-d  \parallel)    
    \leq  W_{n}(s)  \leq  K_{n}(s)  \mbox{,  } \forall s \in G_{n}
    \mbox{.}
   \end{equation}
  \end{lem}	
  \begin{proof}
  For any $x$, $y\in E$, we have
  $$
  W_{n}(x) -W_{n}(y) \leq \sup_{ s \in G_{n} } (\omega(f;d(s,x))-\omega(f;d(s,y))) \mbox{.}
  $$
  By triangular inequality we have $d(s,x) \leq d(s,y)+d(y,x)$.\\
  By growth and subadditivity of $\omega(f)$ we have:\\
  $\omega(f;d(s,x))\leq\omega(f; d(s,y)+d(y,x))\leq\omega(f;d(s,y))+\omega(f;d(x,y))$.\\
  Therefore $W_{n}(x) -W_{n}(y) \leq \omega(f;d(x,y))$.\\
  
  For any $x\in E$ we have\\
  $W_n(x) \leq  \inf_{s\in G_n}K_{n}(f_n)(s)+ \sup_{s\in G_n}\omega(f;d(x,s))$.\\
  By property of moduli of continuity we have\\
  $W_n(x) \leq  \inf_{s\in G_n}K_{n}(f_n)(s)+\sup_{s\in S}f(s)-\inf_{s\in S}f(s)$.\\
  Using (\ref{stable2}) we have $\inf_{s\in G_n}K_{n}(f_n)(s)=\inf_{s\in S_n}f(s)$.\\
  Therefore\\
  $W_n(x) \leq \sup_{s\in S}f(s)+\inf_{s\in S_n}f(s)-\inf_{s\in S}f(s)$\\
  and
  $W_n(x) \leq \sup_{s\in S}f(s)+\omega(f;r_n)$.\\
  Moreover\\
  $W_n(x) \geq  \inf_{s\in G_n}K_{n}(f_n)(s)$.\\
  Using (\ref{stable2}) again we have $\inf_{s\in G_n}K_{n}(f_n)(s)=\inf_{s\in S_n}f(s)$.\\
  Therefore\\
  $W_n(x) \geq \inf_{s\in S}f(s)$. So (\ref{eqmac1.2}) is proved.
  
  Let any $s_0\in G_n$ we have\\
  $W_n(s_0)-K_n(f_n)(s_0) \leq  K_{n}(f_n)(s_0) + \omega(f;d(s_0,s_0))-K_{n}(f_n)(s_0)\leq 0$,\\
  and\\
  $K_n(f_n)(s_0)-W_n(s_0)\leq\sup_{s\in G_{n}}(K_{n}(f_n)(s_0) -K_{n}(f_n)(s) - \omega(f;d(s_0,s)))$.\\
  By $\Omega$-stabilyty of $K_{n}(f_n)$ we have $K_{n}(f_n)(s)-K_{n}(f_n)(s_0))\leq \omega(f;d_n(s_0,s))$\\
  Therefore
  $K_n(f_n)(s_0)-W_n(s_0)\leq \sup_{s\in G_{n}}(\omega(f;d_n(s_0,s))-\omega(f;d(s_0,s)))$\\
  and
  $$
  K_n(f_n)(s_0)-W_n(s_0)\leq \sup_{s\in G_{n}}\omega(f;d_n(s_0,s)-d(s_0,s))
                         \leq \omega(f;\parallel d_n-d\parallel) \mbox{.}
  $$			 
   So (\ref{eqmac1.3}) and  (\ref{eqmac1.4}) are proved.			 
   \end{proof}
  Indeed, from Lemma \ref{lemt.1}, sequence $(W_{n})_{n \in \N}$ is equicontinuous 
  and equibounded. Therefore, by Ascoli's theorem, there exists a subsequence 
  $(W_{\alpha(n)})_{n \in \N}$ which converges to a continuous function denoted by $u$. 
 \begin{thm}\label{th.1}
 The function $u$ is  an $AM\Omega E$ of $f$.
 \end{thm}
 
 \begin{proof}
 We must prove that :\\
 (i) $u$ is an extension of $f$
    \begin{equation}\label{t.6}
  u(s) = f(s) \mbox{ ,} \forall s \in dom(f) \mbox{ ;}
  \end{equation}
  (ii) $u$ is $\Omega-$ optimaly continuous
       \begin{equation}\label{t.7}
  \mid u(x) - u(y) \mid \leq \omega(f;d(x,y)) \mbox{ ,} \forall x,y \in E \mbox{ ;}
  \end{equation}
  (iii) for any open $D \subset E$, such that $D \cap dom(f) = \emptyset $, 
  for any $x \in D$, we have
    \begin{equation}\label{t.8}
     \sup_{y \in \delta D } (u(y) - \omega(u \mid_{\partial D};d(x,y))) 
     \leq u(x) \leq 
     \inf_{y \in \delta D } (u(y) + \omega(u \mid_{\partial D};d(x,y)))    \mbox{.}
  \end{equation} 
For typographical convenience let us assume in the proof that subsequence 
$(W_{\alpha(n)})_{n \in \N}$ is sequence $(W_n)_{n \in \N}$ itself (the true proof 
can easily be restated: replace $n$ by $\alpha(n)$ almost everywhere).

Let us show (i).\\
For any $ \epsilon >0$, there exists $N \in \N$ such that 
$\forall n \geq N$, $ \parallel W_n -u  \parallel_{\infty,E}  \leq \epsilon$.\\
For any $n\geq N$ and $s \in dom(f)$  there exists by (Q1) $s_n\in S_n$ such that $d(s,s_n) \leq r_n$.
 We have\\
 $\mid u(s)-f(s) \mid \leq A_1+A_2+A_3$,\\
 where\\
 $A_1 =\mid u(s)-W_n(s)\mid$, $A_2 =\mid W_n(s)-K_n(s)\mid$,\\
 $A_3 = \mid K_n(s)-K_n(s_n)\mid+\mid f(s_n) -f(s) \mid$.\\
 We have $A_1 \leq \epsilon$. Using (\ref{eqmac1.4}) 
 we have $A_2 \leq \omega(f;\parallel d_n-d\parallel)$.
 Using (\ref{g44.3}) we have \\
 $\mid K_n(s)-K_n(s_n)\mid \leq \omega(f;d_n(s_n,s))
 \leq \omega(f;r_n)+\omega(f;\parallel d_{n}-d  \parallel)$.\\
 In definitive we have\\
  $\mid u(s)-f(s) \mid \leq \epsilon +2\omega(f;r_n)+2\omega(f;\parallel d_{n}-d\parallel)$.\\
 Since this inequality is true $\forall n \geq N$ and $\forall \epsilon >0$ then, using (Q1),(Q3) and
letting $n$ tend to $\infty$, we conclude that $u(s)=f(s)$ so we have proved (i).
 
The proof of (ii) is immediate by letting $n$ tend to $\infty$ in inequality (\ref{eqmac1.1}).

Let us show the right inequality of (iii). Let $D$ an open subset of $E$ such that 
$D \cap dom(f) = \emptyset $. 
For any $\epsilon > 0$, there exists $N \in \N$ such that $ \forall n \geq N$
$ \parallel W_n -u  \parallel_{\infty,E}  \leq \epsilon$.\\
Let $n\geq N$, $x \in D$ and $y \in \partial D$.\\
Using (\ref{Introarlg2.1}) we have\\
$\omega(W_n \mid_{\partial D};d(x,y)) -\omega(u \mid_{\partial D};d(x,y)) 
\leq 
2\parallel W_n -u  \parallel_{\infty,E}$.\\
Now, setting\\
$ A := u(x) - u(y) - \omega(u \mid_{\partial D};d(x,y)) $,\\
we have \\
$ A \leq 4 \epsilon + A_1$,\\
where\\
 $A_1:=W_n(x)-W_n(y)-\omega(W_n \mid_{\partial D};d(x,y))$.\\
Let $D_{n} := \{z \in G_n: b_n(z) \cap D \neq \emptyset \}$  and
 $ \partial D_{n} =\{z \in G_n : b_n(z) \cap \partial D \neq \emptyset \}$. \\
 Since
 $D \subset \cup_{z \in D_{n}}b_n(z)$ and
$\partial D \subset \cup_{z \in \partial D_{n}} b_n(z)$, we have
 $\delta(\partial D,\partial D_n)\leq r_n$ and 
there exists $\tilde{x} \in D_{n}$,$\tilde{y} \in \partial D_{n}$ such that\\
$d(\tilde{x},x) \leq \delta (G_n,E)\leq r_n$ and $d(\tilde{y},y) \leq \delta (G_n,E)\leq r_n$.\\
By lemma \ref{lemt.1}, we have
$$
W_n(x) - K_n(\tilde{x}) \leq   \omega(f;r_n) + 
\omega(f;\parallel d_n -d \parallel) \mbox{ ,}
$$
and
$$
W_n(y) - K_n(\tilde{y}) \leq \omega(f;r_n) + 
\omega(f;\parallel d_n -d \parallel) \mbox{.}
$$
By (\ref{Introarlg2.2}), we have
$$
\omega(W_n\mid_{\partial D};d(x,y)) -\omega(W_n\mid_{\partial D_{n} };d(x,y))
 \leq 4 \omega(W_n ; \delta(\partial D,\partial D_{n})) \mbox{,}
 $$
Since from Lemma \ref{lemt.1}, we have  $\omega(W_n) \leq \omega(f)$, then\\ 
$$
\omega(W_n\mid_{\partial D};d(x,y)) -\omega(W_n\mid_{\partial D_{n} };d(x,y))
 \leq 4\omega(f;\delta(\partial D,\partial D_{n}))\leq 4\omega(f;r_n) \mbox{.}
$$ 
Therefore 
$$
 A_1 \leq A_2+6\omega(f;r_n)+2\omega(f;\parallel d_n -d \parallel) \mbox{.}
$$
where
$$
A_2:=K_n(\tilde{x})-K_n(\tilde{y})-\omega(W_n\mid_{\partial D_{n} };d(x,y)) \mbox{ .}
$$    
Now let us bound $A_2$ from above. We write $A_2:=A_3+A_4$ where \\
$A_3 := 
\omega(K_n\mid_{\partial D_n};d(\tilde{x},\tilde{y}))-
\omega(W_n\mid_{\partial D_{n} };d(x,y))$ and\\
$A_4:= K_n(\tilde{x})-K_n(\tilde{y})-\omega(K_n\mid_{\partial D_n};d(\tilde{x},\tilde{y}))$.\\
We have $d(\tilde{x},\tilde{y})-d(x,y) \leq d(\tilde{x},x)+d(\tilde{y},y) \leq 2r_n$.\\
Using the subadditivity of $\omega(W_n\mid_{\partial D_n})$, we infer that\\ 
$\omega(W_n\mid_{\partial D_n};d(\tilde{x},\tilde{y}))-
\omega(W_n\mid_{\partial D_n};d(x,y)) \leq 2\omega(W_n\mid_{\partial D_n};r_n)
\leq 2\omega(f;r_n)$.\\
Furthermore, using (\ref{Introarlg2.1}) we have\\
$$
\parallel \omega(K_n\mid_{\partial D_n})-\omega(W_n\mid_{\partial D_n})\parallel_{\infty,\R^{+}} 
\leq 2\parallel K_n - W_n \parallel_{\infty,\partial D_n} 
\leq 2\parallel K_n - W_n \parallel_{\infty,G_n} \mbox{.}
$$ 
Now, using (\ref{eqmac1.4}) we have\\
$\parallel K_n - W_n \parallel_{\infty,G_n} \leq \omega(f;\parallel d_n -d \parallel)$.
Therefore
$$
A_3 \leq A_4+2\omega(f;r_n) +\omega(f;\parallel d_n -d \parallel) \mbox{.}
$$

Now, we bound $A_4$ from above.
By theorems \ref{theoExist} and \ref{theoUnicite1}, there exists a unique extension $v$ of 
$K_{n}\mid_{\partial D_n}$ in $G_n$ such that
\begin{equation}\label{problemV.1}
\left\{
\begin{array}{ll}
v(z) =\mu(v;z) &  \forall z \in G_n - \partial D_n \mbox{;}\\
v(z) =K_n(z)     &  \forall z \in \partial D_n \mbox{.}
\end{array}
\right.
\end{equation}
Moreover\\
 $
 v(z) - v(q) \leq \omega(K_{n} \mid_{\partial D_n};d_n(z,q)) \mbox{, }\forall z,q \in G_n\mbox{.}
 $\\
 In particular we have
 \begin{equation}\label{eqAM251}
 v(z) - v(q) - \omega(K_{n} \mid_{\partial D_n};d_n(z,q)) \leq 0 \mbox{, }\forall q \in \partial D_n \mbox{, }
 \forall z \in D_n\mbox{.}
 \end{equation}
 Now, we bound $\sup_{z \in D_n} \mid K_{n}(z)-v(z) \mid$.
 By symmetry we have only to bound $\Delta = \sup_{z \in D_n} (K_{n}(z)-v(z))$ frome above. Let
 $$
 F=\{z \in D_n :K_{n}(z)-v(z)=\Delta \} \mbox{, } M =\sup_{z \in F} K_{n}(z) \mbox{,}
 $$
 and
 $$
 \tilde{F}=\{z \in F:K_{n}(z)=M \} \mbox{.}
 $$
 Let $z_0 \in \tilde{F}$ be such that 
 \begin{equation}\label{eqAM212}
 d(z_0,(S_n \cup \partial D_n)) =\inf_{z\in \tilde{F}}d(z,(S_n \cup \partial D_n)) \mbox{.}
 \end{equation}
 Let us first show that we cannot have $d(z_0,(S_n \cup \partial D_n))>0$. 
  Towards a contradiction let us assume it is the case. We have\\
 $v(z_0) =\mu(v;z_0)$ and $K_n(z_0) =\mu(K_n;z_0)$.  
  Using a similar argument to this of  theorem \ref{theoUnicite}, we infer 
  that $V_n(z_0) \subset \tilde{F}$. Using property (P4) there exists $y \in V_n(z_0)$ 
  such that\\ 
  $d(y,(S_n \cup \partial D_n)) < d(z_0,(S_n \cup \partial D_n))$\\
  which is a contradiction with  definition (\ref{eqAM212}) of $z_0$. \\
   Now, $d(z_0,(S_n \cup \partial D_n))=0$.
   If $z_0 \in \partial D_n$, since $K_n(z_0)=v(z_0)$ we have $\Delta =0$.\\
   If $z_0 \in S_n$ we remark that $z_0 \in D_n$.
   Since $z_0\in S_n \cap D_n$, $S_n\cap D = \emptyset$ and 
   $D_{n} := \{z \in G_n: b_n(z) \cap D \neq \emptyset \}$ 
   there exists $y_0\in D$ such that $d(z_0,y_0) \leq r_n$.\\
   Moreover since $z_0\not\in D$, and by metrical convexity of $E$ 
   there exists $p_0 \in \partial D$ such that 
   $d(z_0,y_0)=d(z_0,p_0)+d(p_0,y_0) \leq r_n$.\\
   By definition of $\partial D_n$ there exists $q_0 \in \partial D_n$ such that 
   $d(p_0,q_0) \leq r_n$.\\
   Therefore
   $d(z_0,q_0)\leq d(z_0,p_0)+d(p_0,q_0)\leq 2r_n$.\\
   We conclude that
   \begin{equation}\label{eqAM22}
   \Delta=
   K_n(z_0)-K_n(q_0)+ v(q_0)-v(z_0)\leq 4\omega(f;r_n) \mbox{.}
    \end{equation}
  From inequalities (\ref{eqAM251}) and (\ref{eqAM22}), we obtain   
 $A_4 \leq  4\omega(f;r_n)$.
 Finally we obtain 
$$
A \leq 4\epsilon + 12\omega(f;r_n)  + 3\omega(f;\parallel d_{n} -d \parallel)\mbox{.}
$$
Since this inequality is true $\forall n \geq N$ then, using (Q1),(Q2),(Q3) and
letting $n$ tend to $\infty$, we conclude that
 $$
 A \leq 4 \epsilon \mbox{ ,}
 $$
which proves the right inequality of (iii). The proof of the left inequality of (iii) is similar but not 
symmetric because of choice of $W_n$. However it leads to similar bounds.
\end{proof}
\section{Numerical tests.}
The tests of this section are done for the following network: $G_n$ is the set of points $(i.h,j.h)$
$i,j =0,...,n$, $h=1/n$ which densifies $\Omega =[0,1]\times [0,1]$ (eventually zoomed and shifted),
 $b_n(x)$ is the ball of center 
$x \in G_n$ radius $h$, $V_n(x)$ the ball of center $x \in G_n$ radius $k.h$. Since norms on $\R^2$ 
are equivalent, balls $b_n(x)$ and $V_n(x)$ can (and will for convenience of implementation), be choosen to be 
those corresponding to $\parallel . \parallel_{\infty}$ or  $\parallel . \parallel_{1}$. Note that, 
for fixed $n$, geodesic metric on $(G_n,V_n)$ will approach euclidean metric on $\Omega$ better 
and better when $k$ increases.
The errors in the following tables are 
$$
e_{n,k}= \sup_{x,y\in G_n}\mid u_{n,k}(x,y)-u(x,y)\mid
$$
where $u_{n,k}$ is the solution of \ref{intro.1} of section 1 for $S=S_n=\partial \Omega \cap G_n$ 
and $f =u\mid S_n$.\\
We first test the algorithm in situations where the solution of the continuous problem is unique and 
known. $u(x,y) =r,\theta,r^{1/2}e^{\theta/2}$ in polar coordinates, $x^{4/3}-y^{4/3}$, in euclidean 
plane.
It is seen that, for a fixed $k$, error becomes stationnary when $n$ increases.\\
Table 7.1: $u(x,y) =r$
\begin{equation}\label{table.00}
\begin{tabular}{|l|r|r|r|r|r|r|r|}
\hline 
k/n    &  8    &   16   &  32     & 64    & 128 & 256  \\
 \hline
1 & 0.023    & 0.023   & 0.023  & 0.023  & 0.023 & 0.02   \\ 
 \hline
2 & 0.0063   & 0.0063  & 0.0066 & 0.0069  & 0.007 & 0.0067   \\ 
 \hline
3 & 0.0062   & 0.0031  & 0.0031 & 0.0031 & 0.0032 & 0.0032    \\ 
\hline
4 & 0.007    & 0.0037  & 0.00205 & 0.0018 &0.0018 & 0.0018    \\ 
 \hline
 5 & 0.0074  & 0.0037  & 0.0021 & 0.001143 &0.00118 &0.0018     \\ 
 \hline
6 & 0.0074   & 0.004   & 0.0022 & 0.001135 & 0.000822 &0.00082    \\ 
 \hline
7 & 0.0079   & 0.004   & 0.0023 & 0.001178 &0.000602& 0.000571 \\ 
 \hline
 \end{tabular}
\end{equation}
Table 7.2: $u(x,y)=\theta$ \\
\begin{equation}\label{table.01}
\begin{tabular}{|l|r|r|r|r|r|r|}
\hline
$k/n$    &  $8$    &   $16$   &  $32$     & $64$    & $128$  &   $256$ \\
 \hline
$1$ & $0.0251$ & $0.0141$  & $0.0139$ &$0.0138$ & $0.0138$ & $0.0138$ \\ 
 \hline
$2$ & $0.165$ & $0.125$  & $0.0347$ &$0.00867$ & $0.00556$ & $0.00421$ \\ 
 \hline
$3$ & $0.236  $ & $0.154 $  & $0.0814$ &$0.0203$ & $0.0054$ & $0.00236$ \\ 
 \hline
 $4$ & $0.244  $ & $0.191 $  & $0.0958$ &$0.0347$ & $0.0088$ & 0.0012 \\ 
 \hline
 \end{tabular}
\end{equation}
Table 7.3: $u(x,y) =r^{1/2}e^{\theta/2}$\\
\begin{equation}\label{table.02}
\begin{tabular}{|l|r|r|r|r|r|r}
\hline
$k/n$    &  $8$    &   $16$   &  $32$     & $64$    & $128$ \\ 
 \hline
$1$ & $0.156$ & $0.112$  & $0.0792$ &$0.0557$ & $0.0402$  \\ 
 \hline
$2$ & $0.22$ & $0.159$  & $0.1123$ &$0.0812$ & $0.0557$  \\
 \hline
$3$ & $0.22  $ & $0.195 $  & $0.1377$ &$0.0971$ & $0.068$  \\ 
 \hline
 \end{tabular}
\end{equation}
Note that we obtain better approximations if we give "thickness $kh$" to the boundary 
that  is if 
we approach the solution of PDE $\Delta_{\infty}u=0$ under Dirichlet's condition 
$u\mid_{\partial_{kh}\Omega} =u_0$ where 
$\partial_{\epsilon}\Omega =[0,1]\times [0,1]-]\epsilon,1-\epsilon[ \times ]\epsilon,1-\epsilon[$.\\
Table 7.4: $u(x,y) =r$\\
\begin{equation}\label{table.03}
\begin{tabular}{|l|r|r|r|r|r|r}
\hline 
k/n  &  8      &   16         &  32     & 64      & 128 \\
 \hline
2 & 0.0037       & 0.0040    & 0.0047  & 0.0057   & 0.0057   \\ 
 \hline
3 & 0.0014       & 0.0015     & 0.0015   & 0.0018 & 0.0023   \\ 
\hline
4 & 0.0000       & 0.0008    & 0.0008 & 0.0008    & 0.00088   \\ 
 \hline
 5 & 0.0000      & 0.0004    & 0.0005 & 0.0005 &  0.0005   \\ 
 \hline
 \hline
 \end{tabular}
\end{equation}

Next we test the algorithm in situations where uniqueness of the solution of the continuous problem 
is  not known: $u_1(x,y) = x^2-y^2$ for $\parallel \mbox{ } \parallel_{\infty}$,\\ 
$u_2(x,y)=\mid x\mid-\mid y\mid$ for  $\parallel \mbox{ } \parallel_{1}$.\\ 
We note that,in these cases, geodesic metric on $G_n$ coincides, for any $k$, with metric on 
$[0,1]\times [0,1]$. So, in these cases, we can take $k=1$. Numerical tests \ref{table.12}  
show that the algorithm computes exactely $u_2$ and that error is linear in $h$ for $u_1$ (see 
Table 7.5).\\
Table 7.5: $u_1(x,y)= x^2 - y^2 $, $u_2(x,y)=\mid x\mid - \mid y\mid$.
\begin{equation}\label{table.12}
\begin{tabular}{|l|r|r|r|r|r|r|}
 \hline
$n$ & $8$    &   $16$   &  $32$     & $64$    & $128$  &   $256$ \\
 \hline
$e1$ & $0.06$ & $0.03$  & $0.015$ &$0.0076$ & $0.0038$ & $0.0019$ \\ 
$e2$ & $0.$ & $0.$  & $0.$ &$0.$ & $0.$ & $0.$ \\ 
\hline
 \end{tabular}
 \end{equation}

To finish we consider the following two examples. In these examples the metrically
convex metric space $(E,d)$ is 
$E=[0,1]\times [0,1]-]1/4-\epsilon,3/4+\epsilon[ \times ]1/2-\epsilon,1/2+\epsilon[$ 
for $\epsilon>0$ small and metric $d$ on $E$ is the geodesic metric constructed from local 
euclidean metric. We set $u_0(x,y)=d((1/2,0),(x,y))$. Figure 7.1 shows function $u_0$ in restriction to 
$G_n$  In the first example we compute the unique solution of 
$\Delta_{\infty}u=0$, $u\mid \Gamma =u_0$ where $\Gamma$  is the union of the boundaries 
of internal and external rectangles. Numerical tests show that the solution in $E$  is 
different from $u_0$: this observation corroborates the fact that geodesic cones are 
not AMLE in general metrically convex metric space (see \cite{ARON2}, appendix).  
\begin{figure}[h]
\caption{$u_0(x,y)$, $n=100$ }
\begin{center}
\fbox{\includegraphics[width=.5\columnwidth, angle=-90]{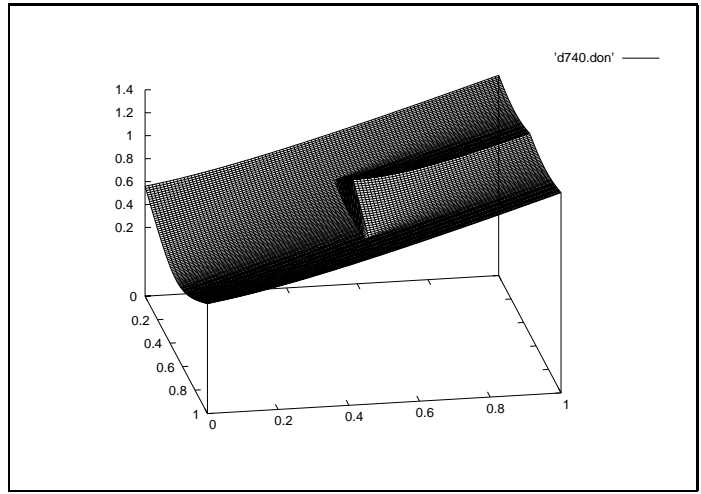}}
\end{center}
\end{figure}

In the second example  we compute  $\Delta_{\infty}u=0$, $u\mid \Gamma =u_0$ where $\Gamma$  
is now the boundary of the external rectangle alone. We obtain a solution 
 which is different 
from $u_0$ and from the solution of the first example. Note the difference between the two examples. 
 In the first one we really compute the solution of PDE $\Delta_{\infty}u=0$ because $E-\Gamma$ is locally euclidean. 
 It is not the case in the second example because the space $(E,d)$ is  not locally euclidean at points 
 of the "free internal boundary". In fact, in this second example it is likely (we are not insured of the 
 convergence of sequence $(W_n)_{n\in \N}$ in Theorem 6.3) that we compute a AMLE of $u_0$. 
 

\section{Appendix.}
As announced in Remark \ref{remar3.7} we prove that
$$
(\frac{1}{2}\sup_{y\in B_h(x)} u(y) +\frac{1}{2}\inf_{y\in B_h(x)} u(y) -u(x))/h^2= -\frac{3}{2}\Delta_{\infty}u(x)+o(h) 
$$
when $u$ is smooth and $Du(x) \neq 0$.\\
 Here $\Delta_{\infty}u =D^2u(D^*u,D^*u)$ and $D^*u=Du/\mid Du\mid$.\\
Since $Du(x) \neq 0$ then, for $h$ sufficiently small, we have $Du(y) \neq 0$ for any $y\in B_h(x)$. Therefore 
$u$ attains its maximum $u^{+}$ and its minimum $u^{-}$ on the boundary of $B_h(x)$. Let us denote $x^{+}$ 
and $x^{-}$ any points of this boundary such that $u(x^{+})=u^{+}$,  $u(x^{-})=u^{-}$. Since $x^{+}$ is a 
maximum of $u(y)$ under the constraint $\mid y-x\mid =h$, vectors $x^{+}-x$ and $Du(x^{+})$ have the same direction 
that is $x^{+}-x=hD^* u(x^{+})$. For the same reason vectors  $x^{-}-x$ and $Du(x^{-})$ have opposite direction that 
is  $x^{-}-x=-hD^* u(x^{-})$.\\
Now, using these expressions of $x^{+}-x$ and $x^{-}-x$ and Taylor formula at $x^+$ and $x^-$ we obtain 
$$
2u(x)=u(x^{+})+u(x^{-})+A+B+h^2o(h)
$$
where $A=h(\mid Du\mid(x^+)-\mid Du\mid(x^-))$ and \\
$B=\frac{1}{2}h^2(D^2u(x^+;D^*(x^+),D^*(x^+))+D^2u(x^-;D^*(x^-),D^*(x^-)))$.\\
Now, using Taylor formula at $x$ for $\mid Du\mid(x^+)$ we have
$$
\mid Du\mid(x^+)= \mid Du\mid(x)+hD(\mid Du\mid)(x;D^*u(x^+))+ho(h) \mbox{.}
$$
By continuity of $y \rightarrow D^*u(y)$, it follows that 
$$
\mid Du\mid(x^+)= \mid Du\mid(x)+hD(\mid Du\mid)(x;D^*u(x))+ho(h) \mbox{.}
$$
Similarly,
$$
\mid Du\mid(x^-)= \mid Du\mid(x)-hD(\mid Du\mid)(x;D^*u(x))+ho(h) \mbox{.}
$$

Since a straightforward computation shows that 
$D(\mid Du\mid)(x;D^*u(x))= \Delta_{\infty}u(x)$,
 and since maps $y \rightarrow D^*u(y)$ and $y \rightarrow D^2u(y)$ are continuous we 
 obtain in definitive
 $$
 2u(x)=u(x^{+})+u(x^{-})+2h^2 \Delta_{\infty}u(x)+h^2\Delta_{\infty}u(x)+h^2o(h) \mbox{,}
 $$ 
which is the announced formula.


 \end{document}